\documentclass[12pt]{amsart}
\usepackage{amsmath, amsthm, amscd, amsfonts}

\newtheorem{theorem}{Theorem}[section]

\newtheorem{proposition}[theorem]{Proposition}
\newtheorem{corollary}[theorem]{Corollary}
\theoremstyle{definition}
\newtheorem{definition}[theorem]{Definition}

\theoremstyle{remark}

\newtheorem{remark}[theorem]{Remark}
\numberwithin{equation}{section}

\newfont{\kh}{msbm10}

\begin{document}
\title[Adjointability of Densely Defined Closed Operators]
{Adjointability of Densely Defined Closed Operators
and the Magajna-Schweizer Theorem}
\author{M. Frank}
\address{Michael Frank, \newline Hochschule f\"{u}r Technik,
Wirtschaft und Kultur (HTWK) Leipzig,
Fachbereich IMN, Gustav-Freytag-Strasse 42A, D-04277 Leipzig,
Germany } \email{mfrank@imn.htwk-leipzig.de}
\author{K. Sharifi}
\address{Kamran Sharifi, \newline Department of Mathematics,
Shahrood University of Technology, P. O. Box 3619995161-316,
Shahrood, Iran} \email{sharifi.kamran@gmail.com and
sharifi@shahroodut.ac.ir}


\begin{abstract}
In this notes unbounded regular operators on Hilbert
$C^*$-modules over arbitrary $C^*$-algebras are discussed. A
densely defined operator $t$ possesses an adjoint
operator if the graph of $t$ is an orthogonal summand. Moreover,
for a densely defined operator $t$ the graph of $t$ is
orthogonally complemented and the range of $P_FP_{G(t)^\bot}$ 
is dense in its biorthogonal complement if and only if $t$ is regular. For a
given $C^*$-algebra $\mathcal A$ any densely defined $\mathcal
A$-linear closed operator $t$ between Hilbert $C^*$-modules is
regular, if and only if any densely defined $\mathcal A$-linear
closed operator $t$ between Hilbert $C^*$-modules admits a
densely defined adjoint operator, if and only if $\mathcal A$ is
a $C^*$-algebra of compact operators. Some further
characterizations of closed and regular modular operators are
obtained.
\end{abstract}
\maketitle

\section{Introduction.}

Hilbert $C^*$-modules are an often used tool in the study of
locally compact quantum groups and their representations, in
noncommutative geometry, in $KK$-theory, and in the study of
completely positive maps between $C^*$-algebras.

A (left) {\it pre-Hilbert $C^*$-module} over a (not necessarily
unital) $C^*$-algebra $\mathcal{A}$ is a left $\mathcal{A}$-module
$E$ equipped with an $\mathcal{A}$-valued inner product $\langle
\cdot , \cdot \rangle : E \times E \to \mathcal{A}$,
which is $\mathcal A$-linear in the first variable and has the
properties:
$$ \langle x,y \rangle=\langle y,x \rangle ^{*},
   \ \ \ \langle x,x \rangle \geq 0 \ \ {\rm with} \
   {\rm equality} \ {\rm if} \ {\rm and} \ {\rm only} \
   {\rm if} \ x=0.
$$
We always suppose that the linear structures of $\mathcal A$
and $E$ are compatible.

A pre-Hilbert $\mathcal{A}$-module $E$ is called a {\it Hilbert
$ \mathcal{A}$-module} if $E$ is a Banach space with respect to
the norm $\| x \|=\|\langle x,x\rangle \| ^{1/2}$. If $E$, $F$
are two Hilbert $ \mathcal{A}$-modules then the set of all
ordered pairs of elements $E \oplus F$ from $E$ and $F$ is a
Hilbert $\mathcal{A}$-module with respect to the $\mathcal A$-valued
inner product $\langle (x_{1},y_{1}),(x_{2},y_{2})\rangle=
\langle x_{1},x_{2}\rangle_{E}+\langle y_{1},y_{2}\rangle _{F}$.
It is called the direct {\it orthogonal sum of $E$ and $F$}.
A Hilbert $\mathcal A$-submodule of a Hilbert $\mathcal A$-module
$F$ is a direct orthogonal summand if $E$ together with its
orthogonal complement $E^\bot$ in $F$ gives rise to an
$\mathcal A$-linear isometric isomorphism of $E \oplus E^\bot$
and $F$. Some interesting results about orthogonally complemented
submodules can be found in \cite{FR2}, \cite{FR3}, \cite{MAG}.
For the basic theory of Hilbert $C^*$-modules we refer to the
book by E.~C.~Lance \cite{LAN} and to respective chapters in the
monographic publications \cite{GVF},  \cite{R-W}, \cite{WEG}.

As a convention, throughout the present paper we assume
$\mathcal{A}$ to be an arbitrary $C^*$-algebra (i.e. not
necessarily unital). Since we deal with bounded and unbounded
operators at the same time we simply denote bounded operators by
capital letters and unbounded operators by lower case letters. We
use the denotations $Dom(.)$, $Ker(.)$ and $Ran(.)$ for domain,
kernel and range of operators, respectively. If $E,\ F$ are
Hilbert $\mathcal{A}$-modules and $W$ is an orthogonal summand in
$E \oplus F$, $P_W$ denotes the orthogonal projection of $E
\oplus F$ onto $W$ and $P_E$ and $P_F$ denote the canonical
projections onto the first and second factors of $E \oplus F$.

Suppose $E,\ F$ are Hilbert $\mathcal{A}$-modules. We denote
the set of all $\mathcal{A}$-linear maps $T: E \to F$
for which there is a map $T^*: F \to E$  such that
the equality
\begin{equation} \label{adjoint}
\langle Tx,y \rangle _{F} = \langle x,T^*y \rangle _{E}
\end{equation}
holds for any $ x \in E,\ y \in F$ by $B(E,F)$.
The operator $T^*$ is called the {\it adjoint operator} of $T$.
The existence of an adjoint operator $T^*$ for some
$\mathcal A$-linear operator $T: E \to F$ implies
that each adjointable operator is necessarily bounded and
$\mathcal{A}$-linear in the sense $T(ax)=aT(x)$ for any $a \in
\mathcal A$, $x \in E$, and similarly for $T^*$. The reason
for this is the requirement that the equality (\ref{adjoint}) is
supposed to hold for any elements of $E$ and $F$, so $T$ has
$E$ as its domain.

In general, bounded $\mathcal{A}$-linear operators may fail
to possess an adjoint operator, however, if $E$ is a full
Hilbert $C^*$-module over a $C^*$-algebra $\mathcal{A}$,
then it is known that each bounded $\mathcal{A}$-linear operator
on $E$ possesses an adjoint operator if and only if $E$ is
{\it orthogonally comparable}, i.e. whenever $E$ appears as a Hilbert
$\mathcal A$-submodule of another Hilbert $\mathcal A$ module $F$
then $E$ is an orthogonal direct summand of $F$ (cf.~\cite{FR2},
Theorem 6.3).

In several contexts where Hilbert $C^*$-modules arise, one also
needs to study 'unbounded adjointable operators', or what are now
known as regular operators. These were first introduced by Baaj
and Julg in \cite{B-J} where they gave an interesting
construction of Kasparov bimodules in $KK$-theory using regular
operators. Later regular operators were reconsidered by
Woronowicz in \cite{WOR}, while investigating noncompact quantum
groups. The functional calculus of regular operators and the
Fuglede-Putnam theorem for Hilbert $C^*$-modules were explained by
Kustermans in \cite{KUS}. Beside these works Kucerovsky gave a
new approach to functional calculus of regular operators in
\cite{KK,KUC}.  Also, Lance gave a brief indication in his book
\cite{LAN} about Hilbert modules and regular operators on them.
Modifying the defining equality (\ref{adjoint}) of adjointability
for unbounded $\mathcal A$-linear operators $t: Dom(t) \subseteq E
\to Ran(t) \subseteq F$ between Hilbert $\mathcal A$-modules $E$
and $F$ the operator $t$ is said to be adjointable if there
exists another $\mathcal A$-linear operator $t^*: Dom(t^*)
\subseteq F \to Ran(t^*) \subseteq E$ such that the equality
\begin{equation}
\langle tx,y \rangle _{F} = \langle x,t^*y \rangle _{E}
\end{equation}
holds for all $ x \in Dom(t),\ y \in Dom(t^*)$. Despite the
good properties of unbounded operators on Hilbert spaces adjointable
unbounded operators on Hilbert $C^*$-modules may lack some good
properties that are wanted in applications. So the notion of
{\it regular operators} was introduced to provide a tractable
class of unbounded $C^*$-linear densely defined closed operators
on Hilbert $C^*$-modules.
An operator $t$ from a Hilbert $\mathcal{A}$-module $E$ to
another Hilbert $\mathcal{A}$-module $F$ is said to be
{\it regular} if
\newcounter{cou001}
\begin{list}{(\roman{cou001})}{\usecounter{cou001}}
\item $t$ is closed and densely defined with domain $Dom(t)$,
\item its adjoint $t^*$ is also densely defined, and
\item the range of $1+t^*t$ is dense in $E$.
\end{list}

Note that as we set $\mathcal{A}= \mathbb{C}$ i.e. if we take
$E$ to be a Hilbert space, then this is exactly the definition of
a densely defined closed operator, except that in that case, both
the second and the third condition follow from the first one. In
\cite{PAL} Pal considered a larger class of operators, semiregular
operators, which are densely defined closable operators whose
adjoints are densely defined. He proved that every closed
semiregular operator (i.e. an operator that satisfies the first
two conditions above) on Hilbert $C^*$-modules over commutative
$C^*$-algebras as well as over subalgebras of $C^*$-algebras of compact
operators is regular (\cite{PAL}, Proposition 4.1, Theorem 5.8).
He also gave an example of a closed semiregular nonregular operator,
and showed that regularity of its adjoint does not ensure regularity
of the original operator (\cite{PAL}, Propositions 2.2 and 2.3).

In the present paper we prove that a densely defined operator $t$
from a Hilbert $ \mathcal{A}$-module $E$ to another Hilbert $
\mathcal{A}$-module $F$ possesses a densely defined adjoint
operator from $F$ to $E$ if the graph of $t$ is orthogonally
complemented in $E \oplus F$ and the range of $ P_F P_{G(t)^{
\perp}}$ is dense in its biorthogonal complement. This fact and
the Magajna-Schweizer theorem show that every densely defined
closed operator on Hilbert $C^*$-modules over $C^*$-algebras of
compact operators is regular, that is for densely defined closed
operators on such Hilbert modules, the second and the third
conditions hold automatically. Magajna, Schweizer and the first
author have presented nice descriptions of $C^*$-algebras of
compact operators in \cite{MAG}, \cite{SCH}, \cite{FR1}. Beside
their work we give further descriptions of such $C^*$-algebras
via some properties of densely defined closed operators.

\section{preliminaries}

In this section we would like to recall some definitions and
present a few simple facts about regular operators on Hilbert
$\mathcal{A}$-modules. For details see chapter 9
and 10 of \cite{LAN}, and the paper \cite{WOR}. We give a
necessary and sufficient condition for closedness of the range
of regular operators.

Let $E, F$  be Hilbert $\mathcal{A}$-modules, we will use the
notation $t : Dom(t) \subseteq E \to F$ to indicate
that $t$ is an $\mathcal{A}$-linear operator whose domain
$Dom(t)$ is a dense submodule of $E$ (not necessarily identical
with $E$) and whose range is in $F$. A densely defined operator
$t: Dom(t) \subseteq E \to F$ is called {\it closed} if its
graph $G(t)=\{(x,tx): \ x \in Dom(t)\}$ is a closed submodule of
the Hilbert $\mathcal{A}$-module $E \oplus F$. In accordance with
the literature we give a stronger definition of
adjointability of densely defined operators that extends the
definition for bounded operators.

\begin{definition} \label{def-strong}
A densely defined operator $t: Dom(t) \subseteq E \to F$
is called {\it adjointable} if it possesses a densely defined map
$t^*: Dom(t^*) \subseteq F \to E$ with the domain
$$
   Dom(t^*)= \{y \in F : {\rm there} \ {\rm exists} \ z \in E \
   {\rm such} \ {\rm that} \ \langle tx,y \rangle _{F} = \langle
   x , z \rangle _{E} \ {\rm for} \ {\rm any} \ x \in Dom(t) \}
$$
which satisfies the property $\langle tx,y \rangle_{F} =
\langle x,t^*y \rangle_{E}$,  for any $x \in Dom(t), \ y
\in Dom(t^*)$.
\end{definition}

The above property implies that $t^*$ is a closed
$\mathcal{A}$-linear map. A densely defined closed
$\mathcal{A}$-linear map $t: Dom(t) \subseteq E \to F$
is called {\it regular} if it is adjointable and the operator
$1+t^*t$ has a dense range. We denote the set of all regular
operators from $E$ to $F$ by $R(E,F)$. There is an alternative
definition of a regular operator between Hilbert $C^*$-modules
(cf.~\cite{WOR}, Definition 1.1), however, Lance has proved in his
book \cite{LAN} that both of them are equivalent. If $t$ is regular
then $t^*$ is regular and $t=t ^{**}$, moreover $t^*t$ is regular and
selfadjoint (cf.~\cite{LAN}, Corollaries 9.4, 9.6 and Proposition
9.9). Define  $Q_{t}=(1+t^*t)^{-1/2}$ and $F_{t}=tQ_{t}$, then
$Ran(Q_{t})=Dom(t)$,  $0 \leq Q_{t} \leq 1$ in $B(E,E)$ and
$F_{t}\in B(E,F)$ (cf.~\cite{LAN}, chapter 9). The bounded
operator $F_{t}$ is called the bounded transform (or $z$-transform)
of the regular operator $t$. The map $t\to F_{t}$ defines a bijection
$$
R(E,F) \to \{ T \in B(E,F):\| T \|\leq 1 \ \, {\rm and} \ \
Ran(1- T^* T ) \ {\rm is} \ {\rm dense \ in} \ F \},
$$
(cf. \cite{LAN}, Theorem 10.4). This map is adjoint-preserving, i.e.
$F_{t}^*=F_{t^*}$, and for the bounded transform
$F_{t}=tQ_{t}=t(1+t^*t)^{-1/2}$ we have $\|F_{t}\|\leq 1$ and
$$
t=F_{t}(1-F_{t}^*F_{t})^{-1/2} \ {\rm and} \
Q_{t}=(1-F_{t}^*F_{t})^{1/2} \, .
$$

Very often there are interesting relationships between regular
operators and their bounded transforms. In fact, for a regular
operator $t$, some properties transfer to its bounded transform
$F_{t}$, and vice versa. Recall the following definitions for a
regular operator $t \in R(E):=R(E,E)$

\begin{itemize}
\item $t$ is called {\it normal} iff $Dom(t)=Dom(t^*)$ and
   $\langle tx,tx \rangle=\langle t^*x,t^*x \rangle $
   for all $x \in Dom(t)$.
\item $t$ is called {\it selfadjoint} iff $t^*=t$.
\item $t$ is called {\it positive} iff $t$ is normal and
   $\langle tx,x \rangle \geq 0$ for all $ x \in Dom(t)$.
\end{itemize}

Then there are the following transfers of properties:
\begin{itemize}
\item $t$ is normal iff $F_{t}$ is normal (cf.~\cite{WOR}, 1.15).
\item $t$ is selfadjoint iff $F_{t}$ is selfadjoint.
\item $t$ is positive iff $F_{t}$ is positive (cf.~\cite{KUS},
      Result 1.14).
\end{itemize}

Let $E$, $F$ be two Hilbert $\mathcal{A}$-modules and suppose that
an operator $T$ in $B(E,F)$ has closed range. We would like to
consider the kernel $Ker(t)$ and the range $Ran(t)$ of $t$.
Closed submodules of Hilbert modules need not to be orthogonally
complemented at all, but Lance states in (\cite{LAN}, Theorem 3.2)
under which conditions closed submodules may be orthogonally
complemented (see also \cite{WEG}, Theorem 15.3.8). For the
special choice of bounded operators $T$ with closed range
one has:
\begin{itemize}
\item $Ker(T)$ is orthogonally complemented in $E$, with
      complement $Ran(T^*)$,
\item $Ran(T)$ is orthogonally complemented in $F$, with
      complement $Ker(T^*)$,
\item the map $T^* \in B(F,E)$ has a closed range, too.
\end{itemize}

The collected facts, as well as Lemmata 4.1 and 4.2 from \cite{N-S}
lead us to the following proposition.

\begin{proposition}
Let $t \in R(E,F)$ and $Ker(t)=\{ x \in Dom(t) \,: tx=0 \}$. Then
\begin{list}{(\roman{cou001})}{\usecounter{cou001}}
\item  $Ker(t)$ and $Ker(t^*)$ both are closed submodules
             of $E$ and $F$, respectively,
\item  $Ran(t)=Ran(F_{t})$ and $Ran(t^*)=Ran(F_{t^*})$,
\item  $Ker(t^*)=Ran(t)^\bot$ and $Ker(t)=Ran(t^*)^\bot$,
\item  $Ker(t)=Ker(F_{t})$ and $Ker(t^*)=Ker(F_{t^*})$.
\item  The regular operator $t$ has closed range if and only if
       its adjoint operator $t^*$ has closed range, and then for
       $|t|:=(t^*t)^{1/2}$ the direct sum decompositions $
                E= Ker(t) \oplus Ran(t^*)= Ker(|t|) \oplus
                \overline{Ran(|t|)}$, $
                F= Ker(t^*) \oplus Ran(t) = Ker(|t^*|) \oplus
                \overline{Ran(|t^*|)}$ hold.
\end{list}
\end{proposition}

\begin{proof}
To show (i), let $\{x_{n}\}$ be a sequence in $Ker(t)$ which
converges to $x \in E$ in norm. Then $t(x_{n})=0$ for any
$n \in \mathbb{N}$. Therefore the sequence $\{t(x_{n})\}$
converges to zero. Closedness of the operator $t$ implies that
$x \in Dom(t)$ and therefore, $tx=0$. So $x\in Ker(t)$ and
$Ker(t)$ is a closed submodule. Since $t$ is regular so
$t^*$ is regular, too, and similarly $Ker(t^*)$ is also a
closed submodule of $F$.

For the proof of (ii) recall that $F_{t}=tQ_{t}$ and $Ran(Q_{t})
=Dom(t)$. Then $Ran(t)=Ran(F_{t})$. Since $t$ is regular
so is $t^*$, thus $Ran(t^*)=Ran(F_{t^*})$.

To demonstrate (iii) we notice that $y \in Ker(t^*)$ if and
only if $\langle tx,y \rangle=\langle x,0 \rangle=0$ for all
$x \in Dom(t)$, or if and only if $y \in Ran(t)^{\bot}$.
Consequently, we have $Ker(t^*)=Ran(t)^{\bot}$. The second
equality follows from the first equality and Corollary 9.4 of
\cite{LAN}.

For the proof of (iv) we know that $Ker(F_{t^*})=
Ran(F_{t})^{\perp}$, cf.~\cite{WEG}, Theorem 15.3.5. Therefore,
$$
  Ker(F_{t^*})=Ran(F_{t})^{\bot}=Ran(t)^{\bot}=Ker(t^*) \, .
$$
Similarly, we obtain $Ker(F_{t})=Ker(t)$.

Finally, we derive (v). The bounded operator $F_{t}$ has closed
range if and only if its adjoint operator $F_{t^*}$ has closed
range. Hence, the regular operator $t$ has closed range if and
only if $t^*$ has. Result 7.19 of \cite{KUS} implies that
$$
\overline{Ran(|t|)}= \overline{Ran(t^{*})}, \ Ker(|t|)=
Ker(t), \  \overline{Ran(|t^*|)}= \overline{Ran(t)}, \
Ker(|t^*|)= Ker(t^*) \, .
$$
The equalities follow immediately from (ii), (iv) and Theorem 3.2
of \cite{LAN}.
\end{proof}

\begin{proposition}
Let $t \in R(E,F)$, then $t$ has closed range if and only if
$Ker(t)$ is orthogonally complemented in $E$ and $t$ is bounded
below on $(Ker(t))^{\perp} \cap Dom(t)$, i.e. $\| tx \| \geq c
\|x \|$,  for all $x \in (Ker(t))^{\perp} \cap Dom(t)$ for a
certain positive constant c.
\end{proposition}

\begin{proof}
Let first $Ran(t)$ be closed then the Proposition 2.2,(v)
implies that $Ker(t)$ is orthogonally complemented in $E$.
We define the $\mathcal{A}$-linear module map
$$
\widetilde{t}: (Ker(t))^{\perp} \cap Dom(t) \to Ran(t)
$$
by $\widetilde{t}x:=tx$ for all $x \in (Ker(t))^{\perp} \cap
Dom(t)$. Then $\widetilde{t}$ is bijection. The inverse of this
mapping exists and is $\mathcal{A}$-linear from $Ran(t)$ into
$(Ker(t))^{\perp}$ with closed domain $Dom(\widetilde{t} \,^{-1})
=Ran(t)$. Moreover $\widetilde{t} \,^{-1}$ is closed since $t$ is.
Therefore, it has to be a bounded operator by the closed graph
theorem, that is there exists a positive constant c such that: $\|
\widetilde{t} \,^{-1} \, x\| \leq c \| x \|$, for each $x \in
Ran\,t$. This implies that $\| tx \| \geq  c^{-1} \|x \|$,  for
all $x \in (Ker\, t)^{ \perp } \cap Dom(t)$.

Conversely, let $t$ be bounded below on $(Ker(t))^{\perp} \cap
Dom(t)$ and $E=Ker(t) \oplus (Ker(t))^{\perp}$. Then if the
sequence $\{y_{n}\} \in Ran(t)$ converges to $y$, there exists a
sequence $\{x_{n}\} \in (Ker(t))^{\perp} \cap Dom(t)$ such that
$y_{n}=t(x_{n})$. Then $(x_{n}-x_{m}) \in (Ker(t))^{\perp}$, and
therefore $\|x_{n}-x_{m}\| \leq c^{-1} \|y_{n}-y_{m}\|$ converges
to zero as $m,n$ go to infinity. This means that there exists an
element $x\in (Ker(t))^{\perp}$ such that the sequence $\{ x_{n}
\}$ converges to $x$ in norm and the sequence $\{ t (x_{n}) \}$
converges to $y$ in norm. The closedness of $t$ implies that $x
\in Dom(t)$ and $tx=y$.
\end{proof}

\section{Adjointability of densely defined operators}

Corollary 2.4 of \cite{FR1} shows that a bounded
$\mathcal{A}$-linear operator $T:E \to F$ possesses an adjoint
operator $T^*:F \to E$ if and only if the graph of $T$ is an
orthogonal summand of the Hilbert $\mathcal{A}$-module $E \oplus
F$. This fact motivates us to give a sufficient condition for the
adjointability of densely defined operators via their graphs.

\begin{theorem}
Let $E,F$ be two Hilbert $\mathcal{A}$-modules and $t:
Dom(t)\subseteq E \to F$ be a densely defined operator. If the
graph of $t$ is orthogonally complemented in $E \oplus F$ and the
range of $ P_F P_{G(t)^{ \perp}}$ is dense in its biorthogonal
complement then $t$ is adjointable. In this case $t$ is closed
and $1+t^*t$ is surjective.
\end{theorem}

\begin{proof}
Consider the unitary element $V$ of $B(E \oplus F,F\oplus E)$
defined by $V(x,y)=(y,-x)$. Then $G(t)$ is orthogonally
complemented in $E \oplus F$ if and only if
$V(G(t))$ is orthogonally complemented in $F\oplus E$.
The orthogonal complement of $V(G(t))$ with respect to
$F\oplus E$ is the closed set
$$
[V(G(t))]^{\perp}= \{(y,z): y \in F, \ z \in E, \ {\rm such}
\ {\rm that} \ \langle(tx,-x),(y,z)\rangle=0, \ {\rm for}
\ {\rm all} \ x \in Dom(t) \}
$$
$$
\ \ \ \ \ = \{(y,z):y \in F, \ z \in E, \ {\rm such}
\ {\rm that} \ \langle tx,y\rangle=\langle x,z \rangle,
\ {\rm for} \ {\rm all} \ x \in Dom(t) \} \, .
$$
Now we define
$$
Dom(t^*):= \{y \in F : {\rm there} \ {\rm exists} \
z \in E \ {\rm such} \ {\rm that} \ \langle tx,y \rangle  =
\langle x,z \rangle \ {\rm for} \ {\rm all} \ x \in Dom(t) \} \, .
$$
The set $Dom(t^*)$ is a non-trivial submodule of $F$ since the set
$[V(G(t))]^{\perp}$ is a non-trivial submodule of $F \oplus E$.
The domain of $t$ is dense in $E$, and so we consider elements
$y \in F$ such that an element $z$ with the property
$$
\langle tx,y \rangle  = \langle x,z \rangle, \ {\rm for}
\ {\rm any} \ x \in Dom(t)
$$
exists and is unique. This set is not empty since $y=0$ and $z=0$
forms an admissible pair of elements. Collecting all such elements
$y \in F$ we can define an operator $t^*: Dom(t^*)\subseteq F \to
E$ by $t^*\, y=z$. Clearly $t^*$ is $ \mathcal{A}$-linear and
satisfies $\langle tx,y \rangle  = \langle x , t^*y \rangle$ for
all $x \in Dom(t), \ y \in Dom(t^*)$. Moreover we have
$[V(G(t))]^{\perp}= \{(y,t^*y): \ y \in Dom(t^*) \}=G(t^*)$, i.e.
$F \oplus E=V(G(t)) \oplus G(t^*)$. Since $Ran( P_F P_{G(t)^{
\perp}} ) =Dom(t^*)$ and $( Ran( P_F P_{G(t)^{ \perp}} ) )^{
\perp}= Ker( P_{G(t)^{ \perp}} P_{F}^{*})=\{0 \}$, we find $
\overline{Dom(t^*)} = \{ 0 \} ^{ \, \perp}=F$. The graph of $t$
is orthogonally complemented, so Lemma 15.3.4 of \cite{WEG}
implies that $G(t)$ is closed in $E \oplus F$. Suppose $u \in E$
is an arbitrary element then $(0,u) \in F \oplus E = V(G(t))
\oplus G(t^*)$, that is, there exist elements $ x_0 \in Dom(t)$
and $y_0 \in Dom(t^*)$ such that $ y_0=-t(x_0 )$ and $
u=t^*(y_0)- x_0$, consequently, $x_0 \in Dom(t^*t)$ and
$u=-(1+t^*t)x_0 $. Hence, $1+t^*t$ is surjective.
\end{proof}

The above theorem shows that the closedness and the assumption to
unbounded $C^*$-linear densely defined operators of possessing a
densely defined adjoint operator can be reduced in some results of
\cite{KK, LAN}. On the contrary to the situation for bounded
operators, the converse of the above theorem is not valid. An
example of a selfadjoint densely defined closed operator which
graph is not orthogonally complemented was given by Hilsum in
\cite{HIL} (see also \cite{LAN}, page 103). Now we can find a
necessary and sufficient condition as follows:

\begin{corollary}
Let $t : Dom(t) \subseteq E \to F$ be an $\mathcal A$-linear
densely defined operator between Hilbert $\mathcal A$-modules $E$
and $F$. Then the graph of $t$ is orthogonally complemented in $E
\oplus F$ and $\overline{Ran ( P_F P_{G(t)^{ \perp}})}= Ran ( P_F
P_{G(t)^{ \perp}}) ^{ \, \perp \perp}$ if and only if $t$ is
regular, i.e. $t$ is adjointable, closed, and the range of ${\rm
1}+t^*t$ is dense in $E$.
\end{corollary}

\begin{proof}The assertion is a direct conclusion from Theorem 9.3
of \cite{LAN} and Theorem 3.1 above.
\end{proof}

The criterion found by Kucerovsky as Proposition 6 in \cite{KK}
now reads as follows:

\begin{corollary}
Let $t: Dom(t) \subseteq E \to F$ be an $\mathcal A$-linear
densely defined operator between Hilbert $\mathcal A$-modules $E$
and $F$. The operator $t$ is regular if and only if $t$ is
adjointable, closed, and for any positive real number $c$ the
operator $c{\rm 1} + t^*t$ is bijective. If $t$ is selfadjoint,
then $t$ is necessarily closed and so $t$ is regular if and only
if the operator $ci \pm t$ is bijective for any non-zero real
constant $c$.
\end{corollary}

In \cite{KUC} Kucerovsky has given a geometrical criterion for
regularity of closed operators,  cf. (\cite{KUC}, Proposition 5).
Here, by reducing some of his suppositions we can sharpen his
criterion.

\begin{theorem}
Let $E$, $F$ be two Hilbert $\mathcal{A}$-modules. Then a closed
operator is regular if and only if there exists a Hilbert
$\mathcal{A}$-module $G$ and a bounded adjointable operator $S
\in B(G,E \oplus F)$ such that
\begin{list}{(\roman{cou001})}{\usecounter{cou001}}
\item   the graph of the operator is the range of $S$,
\item   the range of $P_{E}S$ is dense $E$, and
\item   the range of $P_{F}P_{Ker(S^*)}$ is dense in its biorthogonal complement.
\end{list}
\end{theorem}

\begin{proof}
Let a regular operator $t: Dom(t)\subseteq E \to F$ be given.
Suppose $G$ is the graph of $t$ and $S$ is the inclusion of $G$
into $E \oplus F$. The graph of $t$ is orthogonally complemented
in $E \oplus F$ and so $S$ is adjointable and $Ker(S^*)=Ran(S) ^{
\perp}= G(t)^{ \perp}$. Furthermore, $Ran(P_{E}S)=Dom(t)$ and
$Ran(P_F P_{G(t)^{ \perp}})=Dom(t^*)$ are dense in $E$ and $F$,
respectively.

Conversely, let $ t: Dom(t)\subseteq E \to F$ be closed and
suppose the conditions (i), (ii) and (iii) hold. Then the range of
$S$ is closed, hence Theorem 3.2 of \cite{LAN} implies that it is
an orthogonal summand, with complement $Ker(S^*)$. The range of
$P_{E}S$ is dense in $E$, so $t$ is a densely defined closed
operator whose graph is orthogonally complemented in $E \oplus F$
and $\overline{Ran ( P_F P_{G(t)^{ \perp}})}= Ran ( P_F P_{G(t)^{
\perp}}) ^{ \, \perp \perp}$, that is, $t$ is regular by Corollary
3.2.
\end{proof}

Suppose that $\mathcal{A}$ is an arbitrary $C^*$-algebra of
compact operators. It is well-known that $\mathcal A$ has to be
of the form $\mathcal{A}$=$c_{0}$-$ \oplus_{i \in I}\mathcal{K}
(H_{i})$, i.e. $\mathcal{A}$ is a $c_{0}$-direct sum of elementary
$C^*$-algebras $\mathcal{K}(H_{i})$ of all compact operators
acting on Hilbert spaces $H_{i}, \ i \in I$ (cf.~\cite{ARV},
Theorem 1.4.5).

Magajna and Schweizer have shown, respectively, that
$C^*$-algebras of compact operators can be characterized by the
property that every norm closed (coinciding with its
biorthogonal complement, respectively) submodule of every Hilbert
$C^*$-module over them is automatically an orthogonal summand, cf.
\cite{MAG}, \cite{SCH}. Recently further generic properties of the
category of Hilbert $C^*$-modules over $C^*$-algebras which
characterize precisely the $C^*$-algebras of compact operators
have been found by the first author in \cite{FR1}. We recall
results by Magajna, Schweizer and Frank as
follows:

\begin{theorem}
Let $\mathcal{A}$ be a $C^*$-algebra. The
following conditions are equivalent:
\begin{list}{(\roman{cou001})}{\usecounter{cou001}}
\item $\mathcal{A}$ is an arbitrary $C^*$-algebra of compact operators.

\item For every Hilbert $\mathcal{A}$-module $E$ every
Hilbert $\mathcal{A}$-submodule $F \subseteq E$ is automatically
orthogonally complemented, i.e. $F$ is an orthogonal summand.

\item For every Hilbert $\mathcal{A}$-module $E$ Hilbert
$\mathcal{A}$-submodule $F \subseteq E$ that coincides
with its biorthogonal complement $F^{\perp \perp}
\subseteq E$ is automatically orthogonally complemented in $E$.

\item For every pair of Hilbert $\mathcal{A}$-modules
$E, F$, every bounded $\mathcal{A}$-linear map $T:E
\to F$ possesses an adjoint bounded $\mathcal{A}$-linear map
$T^*:F \to E$.

\item The kernels of all bounded $ \mathcal{A}$-linear
operators between arbitrary Hilbert $ \mathcal{A}$-modules are
orthogonal summands.

\item The images of all bounded $ \mathcal{A}$-linear
operators with norm closed range between arbitrary
Hilbert $\mathcal{A}$-modules are orthogonal summands.

\item For every Hilbert $ \mathcal{A}$-module $E$
every Hilbert $ \mathcal{A}$-submodule is automatically
topologi\-cally complemented there, i.e. it is a topological
direct summand.

\item For every (maximal) norm closed left ideal $I$
of $ \mathcal{A}$ the corresponding open projection $p
\in \mathcal{A}^{**}$ is an element of the multiplier
$C^*$-algebra $M(\mathcal{A})$ of $ \mathcal{A}$.
\end{list}
\end{theorem}

Consider the $C^*$-algebra of compact operators as a Hilbert
$C^*$-module over itself. Pal has proved in Theorem 5.8 of
\cite{PAL} that every closed semiregular operator (i.e. every
densely defined closed operator which adjoint is densely defined)
on Hilbert $C^*$-modules over $C^*$-algebras of compact operators
is regular. Corollary 3.2 and the second part of the above theorem
give a short proof for Pal's Theorem. Moreover we can reformulate
Pal's Theorem as follows:

\begin{remark}
In the category of all Hilbert $C^*$-modules over a $C^*$-algebra
of compact operators every densely defined closed $C^*$-linear
operator between Hilbert $C^*$-modules is regular.
\end{remark}

\begin{corollary}Let $\mathcal{A}$ be a $C^*$-algebra. The
following conditions are equivalent:
\begin{list}{(\roman{cou001})}{\usecounter{cou001}}
\item $\mathcal{A}$ is an arbitrary $C^*$-algebra of
compact operators.
\end{list}
\begin{list}{(\roman{cou001})}{\usecounter{cou001}}
   \setcounter{cou001}{8}
\item For every pair of Hilbert $ \mathcal{A}$-modules
$E, F$, every densely defined closed operator $t:
Dom(t)\subseteq E \to F$ possesses a densely defined
adjoint operator $t^*: Dom(t^*)\subseteq  F \to E$.

\item For every pair of Hilbert $ \mathcal{A}$-modules
$E, F$, every densely defined closed operator $t:
Dom(t)\subseteq E \to F$ is regular.

\item The kernels of all densely defined closed
operators between arbitrary Hilbert $\mathcal{A}$-modules
are orthogonal summands.

\item The images of all densely defined closed operator
with norm closed range between ar\-bitrary Hilbert
$\mathcal{A}$-modules are orthogonal summands.
\end{list}
\end{corollary}

\begin{proof}
Theorem 3.1, Pal's Theorem and condition (ii) imply (ix), (x),
(xi) and (xii). To show the contrary let condition (ix) hold and
let $T: E \to F$ be an arbitrary bounded $\mathcal{A}$-linear map
between Hilbert $\mathcal{A}$-modules $E$ and $F$. The operator
$T: Dom(T)=E \to F$ is a densely defined closed operator (since it
is bounded), and so condition (ix) implies that there exists a
(possibly unbounded) densely defined operator $T^*:
Dom(T^*)\subseteq F \to E$ such that $$ \langle Tx,y \rangle =
\langle x,T^*y \rangle, \ \ {\rm for} \ {\rm all} \ x \in
Dom(T)=E, \ y \in Dom(T^*) \, . $$ Then $T^*$ is bounded on the
pre-Hilbert module $Dom(T^*)$. The domain of $T^*$ is dense in $F$
and $E$ is a Hilbert module, so $T^*$ has a unique bounded
$\mathcal{A}$-linear extension $\widetilde{T^{*}}: F \to E$ such
that  $\langle \, Tx,y \rangle = \langle x, \widetilde{T^{*}}\,y
\rangle$ for all $ x \in E, \ y \in F$. Therefore every bounded
$\mathcal{A}$-linear map $T:E \to F$ possesses an adjoint bounded
$\mathcal{A}$-linear map, i.e. condition (iv) holds. Condition (x)
implies (ix) and hence, (iv). Conditions (xi) and (xii) imply
conditions (v) and (vi), respectively, since each everywhere
defined bounded operator is a densely defined closed operator.
\end{proof}

\begin{corollary}
Suppose $\mathcal{A}$ is any $C^*$-algebra which does not admit a
faithful $*$-representation as a $C^*$-subalgebra in some
$C^*$-algebra of compact operators. Then there exists a densely
defined closed operator $t$ between two full Hilbert $C^*$-modules
over $\mathcal{A}$ such that $t$ is not regular or, even more,
the adjoint operator $t^*$ of $t$ is not densely defined (and therefore,
there does not exist any adjoint operator in the strong sense of
Definition \ref{def-strong}).
\end{corollary}

This fact directly follows from Corollary 3.7,(i),(ix) and (x). In
other words, Hilsum's example reflects a quite regular case for
unbounded densely defined closed operators between Hilbert
$C^*$-modules over non-compact $C^*$-algebras of any kind.

Let $ \mathcal{K}(H)$ be the $C^*$-algebra of all compact
operators on a Hilbert space $H$. Let $e \in \mathcal{K}(H)$ be
an arbitrary minimal projection and $E$ be a $
\mathcal{K}(H)$-module. Suppose $E_{e}:=e E=\{ex: \ x\in E \}$,
then $E_{e}$ is a Hilbert space with respect to the inner product
$(.,.)=trace\,(\langle.,.\rangle)$, which is introduced in
\cite{B-G}. Let $B(E)$ and $B(E_{e})$ be $C^*$-algebras of all
bounded adjointable operators on Hilbert $ \mathcal{K}(H)$-module
$E$ and Hilbert space $E_{e}$, respectively. Baki\'c and
Gulja\v{s} have shown that the map $ \Phi : B(E) \to
B(E_{e}) , \ \Phi(T)= T| _{E_{e}} $ is a $*$-isomorphism of
$C^*$-algebras (cf. \cite{B-G}, Theorem 5). On the other hand,
Remark 3.4 and Theorem 10.4 of \cite{LAN} give adjoint-preserving
bijection maps $t \to F_{t}=t(1+t^*t)^{-1/2}$ as
follows:
$$
R(E) \ \, \to \ \{ T \in B(E):\| T \|\leq 1 \ \, {\rm and} \ \
Ran(1- T^*T ) \ {\rm is} \ {\rm dense} \ in \ E \} \, ,
$$
$$
R(E_{e}) \to \{ T \in B(E_{e}):\| T \|\leq 1 \ \, {\rm and} \ \
Ran(1- T^*T ) \ {\rm is} \ {\rm dense} \ {\rm in} \ E_{e} \} \, .
$$

These maps together with the $*$-isomorphism $\Phi$ give an
adjoint-preserving bijection map between all densely defined
closed operators on Hilbert $\mathcal{K}(H)$-module $E$ and all
densely defined closed operators on Hilbert space $E_{e}$. It
means that all densely defined closed operators on a Hilbert
$\mathcal{K}(H)$-module $E$, are reduced by a suitable Hilbert
space contained in $E$.

{\bf Acknowledgement}: The authors would like to thank the referee
for his/her useful comments.

{\bf Added in Proof:} After the present paper has been accepted we
learned about a related publication \cite{Gu2} by B. Gulja\v{s}.
His results are rather complementary and describe situations for
C*-algebras of compact operators. The authors thank 
B.~Gulja{\v{s}} for pointing out a gap in the initial proof of
Theorem 3.1.

\end{document}